\documentclass[12pt]{article}
\usepackage{amssymb,amsmath,amsthm,epsfig}
\usepackage{latexsym, enumerate}
\usepackage{eepic}
\usepackage{epic}
\usepackage{graphicx}
\usepackage{color}
\usepackage{ifpdf}
\usepackage{subfigure}
\usepackage{tikz}
\usepackage{dsfont}
\usepackage{multirow}
\usepackage{makecell}
\usepackage{algorithm}
\usepackage{bm}

\theoremstyle{plain}

\theoremstyle{definition}

\theoremstyle{remark}

\begin{document}
\title{ \large\bf Uniqueness of the inverse reaction coefficient problems for nonlocal diffusion models}

\author{
Guang-Hui Zheng\thanks{Corresponding author. College of Mathematics and Econometrics, Hunan University, Changsha 410082, Hunan Province, China. Email: zhenggh2012@hnu.edu.cn}
\and
Ming-Hui Ding\thanks{College of Mathematics and Econometrics, Hunan University, Changsha 410082, Hunan Province, China. Email: minghuiding@hnu.edu.cn}
}

\date{}
\maketitle

\begin{center}{\bf ABSTRACT}
\end{center}\smallskip
In this paper, we consider the inverse reaction coefficient problems (IRCPs) for nonlocal diffusion equation and multi-term time-fractional nonlocal diffusion equation from the average nonlocal flux data in external reaction region. Based on the nonlocal maximum principle we established, the uniqueness theorem for IRCPs are proved.

\smallskip
{\bf keywords:} Nonlocal diffusion, reaction coefficient, nonlocal maximum principle, uniqueness, average nonlocal flux

\section{Introduction}
Nonlocal models and nonlocal diffusion operators are widely applied in many fields, such as continuum mechanics \cite{Marta D'Elia2013, S. Silling2000}, biology \cite{Carrillo2005, Mogilner1999}, jump process \cite{M.T. Barlow2009,R.F. Bass2010,N. Burch2011}, graph theory \cite{L. Lovasz2006}, image analyses \cite{A. Buades.2010,G. Gilboa2007,Y. Lou2010}, machine learning \cite{L. Rosasco2010}, and phase transitions \cite{P. Bates1999,P. Fife2010}.

The difference between the nonlocal model and the classical partial differential equation model is that in the latter case, the interaction between two regions occurs only because of contact, while in the former case, the interaction can occur at a certain distance. Let $\Omega$ be a bounded domain in $\mathbb{R}^d$ $(d\geq1)$, and $u(\bm{\mathrm{x}}):\Omega\rightarrow \mathbb{R}$ define the action of the nonlocal diffusion operator $\mathcal{L}$ on the function $u(\bm{\mathrm{x}})$ as follows\\
\begin{eqnarray}
\label{mi1}
\begin{split}
\mathcal{L}u(\bm{\mathrm{x}}):=2\int_{\mathbb{R}^d}(u(\bm{\mathrm{y}})-u(\bm{\mathrm{x}}))\gamma (\bm{\mathrm{x}},\bm{\mathrm{y}})d{\bm{\mathrm{y}}},\ \ \ \forall \bm{\mathrm{x}}\in \Omega\subseteq \mathbb{R}^d,
\end{split}
\end{eqnarray}
here, the kernel $\gamma(\bm{\mathrm{x}},\bm{\mathrm{y}}):\Omega\times\Omega\rightarrow \mathbb{R}$ is a non-negative symmetric function, and satisfies the following inequalities
\begin{eqnarray}
\label{kernel function}
\begin{split}
\gamma^*\geq\gamma(\bm{\mathrm{x}},\bm{\mathrm{y}})\mid \bm{\mathrm{x}}-\bm{\mathrm{y}}\mid^{d+2\beta}\geq\gamma_{*},\ \ \ \text{for}\ \ \bm{\mathrm{y}}\in B_\varepsilon(\bm{\mathrm{x}}),
\end{split}
\end{eqnarray}
and
\begin{eqnarray}
\label{kernel function1}
\begin{split}
\gamma(\bm{\mathrm{x}},\bm{\mathrm{y}})=0,\ \ \ \text{for}\ \ \bm{\mathrm{y}}\in \mathbb{R}^d\setminus B_\varepsilon(\bm{\mathrm{x}}),
\end{split}
\end{eqnarray}
where $B_\varepsilon(\bm{\mathrm{x}}):=\{\bm{\mathrm{y}}\in \mathbb{R}^d: |\bm{\mathrm{y}}-\bm{\mathrm{x}}|<\varepsilon\}$, $\beta\in (0,1)$, $\gamma^*$ and $\gamma_*$ is positive constants. For nonlocal operator $\mathcal{L}$, the value of $\mathcal{L}u$ at $\bm{\mathrm{x}}$, all information about $\bm{\mathrm{y}}\neq \bm{\mathrm{x}}$ is required, and the value of $\Delta u$  at $\bm{\mathrm{x}}$ which only needs  information at $\bm{\mathrm{x}}$ for local operators (see \cite{Du2012}).

Next we consider the operator $\mathcal{L}$ is due to its participation in nonlocal diffusion equation (NDE) with Dirichlet volume-constrained problem
\begin{eqnarray}
\label{mix1}
\begin{cases}
\begin{split}
\frac{\partial u}{\partial t}-\mathcal{L}u+q(\bm{\mathrm{x}})u&=\phi(\bm{\mathrm{x}})v(t),\ \ \ \text{in}\ \ \Omega\times(0, T),\\
u&=0,\ \ \ \text{on}\ \ (\mathbb{{R}}^d\setminus\Omega)\times(0, T),\\
u&=0,\ \ \ \text{on}\ \ \Omega\times\{0\},\\
\end{split}
\end{cases}
\end{eqnarray}
where the reaction coefficient $q\in C(\overline{\Omega})$, and $q\geq 0$, the input source is formed by the separated variables $\phi(\bm{\mathrm{x}})v(t)$, where $v(t)$ is the time-varying strength of source, and $\phi(\bm{\mathrm{x}})$ represents the space-position information. The Dirichlet volume constraints are natural extensions, to the nonlocal case, of Dirichlet boundary condition for classical diffusion problem.

Since the time-fractional diffusion equation is closely related to fractional Brownian motion, and is an important tool for describing anomalous diffusion in highly heterogeneous media \cite{I. Podlubny1999, A.A.Kilbas2006}. We also consider the following multi-term time-fractional nonlocal diffusion equation (MTTFNDE)
\begin{eqnarray}
\label{mix2}
\begin{cases}
\begin{split}
{_{0}D}_t^{\alpha}u+\sum_{k=1}^{m}p_{k}({_{0}D}_t^{\alpha_k}u)-\mathcal{L}u+q(\bm{\mathrm{x}})u&=\phi(\bm{\mathrm{x}})v(t),\ \ \ \text{in}\ \ \Omega\times(0, T),\\
u&=0,\ \ \ \text{on}\ \ (\mathbb{R}^d\setminus\Omega)\times(0, T),\\
u&=0,\ \ \ \text{on}\ \ \Omega\times\{0\},\\
\end{split}
\end{cases}
\end{eqnarray}
where $m$ is a fixed positive integer, and $p_k$ are positive constants. The fractional orders satisfy $0<\alpha_1<\alpha_2<\cdots<\alpha_m<\alpha<1$, and $_{0}D_t^{\alpha}u$ is the Caputo fractional derivative defined by \cite{I. Podlubny1999}
\begin{eqnarray*}
\begin{split}
_{0}D_t^{\alpha}u=\frac{1}{\Gamma(1-\alpha)}\int_0^t(t-s)^{-\alpha}\frac{\partial u(\bm{\mathrm{x}},s)}{\partial s}ds,
\end{split}
\end{eqnarray*}
and $\Gamma(\cdot)$ denotes the Gamma function.

As for the direct problems for nonlocal diffusion models, i.e., the volume-constrained problem, which have been studied extensively in the past few years \cite{Du2019,X. Tian2013,H. Tian2015,R. Nochetto2015,A. Bueno-Orovio2014,Q. Yang2011}. However, about the corresponding inverse problems, the results are very limited (see \cite{S. Tatar2016,J. Jia2018,M. D'Elia2013}). In this paper, our goal is to identify the reaction coefficient $q(\bm{\mathrm{x}})$ for NDE and MTTFNDE from the average nonlocal flux measurement data, which are usually measured on the accessible part of the external interaction region in nonlocal models. The solutions to system (\ref{mix1}) and (\ref{mix2}) will be denoted by $u_j^i(\bm{\mathrm{x}},t;q)$, in order to indicate its dependence on the reaction coefficient $q$, and correspond to the input sources $\phi_j(\bm{\mathrm{x}})v_i(t)$, $i=1,2$; $j=1,2,\cdots$. Hereafter, $C_0^{2,1}(\Omega\times(0,T))$ denotes the function space in which the functions are 2-times continuously differentiable with respect to spatial variable and 1-times continuously differentiable with respect to time variable, and vanish near the boundary of $\Omega\times(0,T)$. $C$ refers to a generic constant which may differ at different occurrences.

\emph{Inverse reaction coefficient problem (IRCP) for NDE and MTTFNDE}: Given the input source $\phi_j(\bm{\mathrm{x}})v_i(t)$, $i=1,2$; $j=1,2,\cdots$, the average nonlocal flux data set are
\begin{eqnarray}
\small
\label{me1}
\begin{split}
\int_{\Omega_a\times(0,T)}\mathcal{N}(\bm{\Theta}\cdot\mathcal{D}^* (u_j^i(\bm{\mathrm{x}},t;q)))h(s,t)dsdt,
\end{split}
\end{eqnarray}
determine the reaction coefficient $q$ (see Figure 1.1 for a schematic illustration). Here $\Omega_a\subseteq (\mathbb{R}^d\setminus \Omega)$ is an accessible region of $\mathbb{R}^d\setminus \Omega$, the operator $\mathcal{N}$ is the nonlocal interaction operator, $\bm{\Theta}$ is the second-order symmetric positive definite tensor, $\mathcal{D}^*$ is the adjoint operator of nonlocal divergence operator $\mathcal{D}$ (see Section 2). $h$ is a nonzero nonnegative function, which can be interpreted as a characterization of measure instrument.\\

\begin{figure}[h]
\label{fig}
\begin{center}
\begin{tikzpicture}[line width=1.5pt,scale=0.7]
\draw[very thick] (0,0) to [out=110,in=185] (2,2) to [out=15,in=140] (4,2) to [out=-30,in=50] (5,-1) to [out=-150,in=-50] (0,0);
\draw[very thick] (7,0) to [out=90,in=190] (8,1) to [out=190,in=140] (9,1) to [out=-30,in=50] (10,-1) to [out=-150,in=-80] (7,0);
\node at (4,1.5){$\mathbf{\Omega}$};
\node at (6.3,2){$\mathds{R}^d\setminus\mathbf{(\Omega\cup\Omega_a)}$};
\node at (7.7,0.2){$\mathbf{\Omega_{a}}$};
\node at (1.5,0.5){$\bm{\otimes}$};
\node at (2,-0.5){$\bm{\otimes}$};
\node at (4,0.5){$\bm{\otimes}$};
\node at (2.5,1){$\bm{\otimes}$};
\node at (3.5,-0.5){$\bm{\otimes}$};
\node at (8.2,-0.5){$\bm{\odot}$};
\node at (9.3,0){$\bm{\odot}$};
\node at (8.6,0.5){$\bm{\odot}$};
\node at (9,-0.8){$\bm{\odot}$};
\end{tikzpicture}
\end{center}
\caption{The schematic of IRCP, $\Omega$ is the physical domain, $\Omega_{a}$ denotes the accessible region, and $\bm{\otimes}$ represent the input source locations, $\bm{\odot}$ indicate the measurement locations.}
\end{figure}
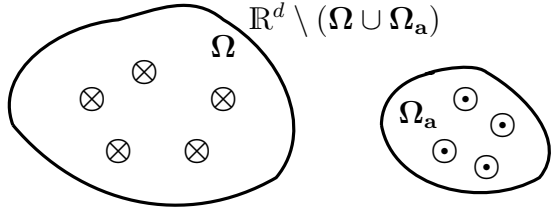

The main results of this paper read as follows.\\

\textbf{Theorem 1.} Let $\{\phi_j\}_{j=1}^{\infty}\in C(\Omega)$ be a complete set in $L^2(\Omega)$, $v\in C^1(0,T)$ and $h\in C_0^{2,1}(\Omega_a\times(0,T))$ be given nonzero nonnegative functions, and $v$ satisfies $v(0)=0$. Assume $p(\bm{\mathrm{x}})$, $q(\bm{\mathrm{x}})\in C(\overline\Omega)$, $p$, $q\geq 0$ on $\Omega$. Let $u_j^i(\bm{\mathrm{x}},t;p)$, $u_j^i(\bm{\mathrm{x}},t;q)$ be the classical solutions of problem (\ref{mix1}) corresponding to the input sources $\phi_j(\bm{\mathrm{x}})v_i(t)$ ($i=1,2$; $j=1,2,\cdots$) with the reation coefficients $p$ and $q$ respectively. If we choose $v_1=v$, $v_2=v'(t)$ such that
\begin{eqnarray}
\small
\label{condtion}
\begin{split}
\int_{\Omega_a\times(0,T)}\mathcal{N}(\bm{\Theta}\cdot \mathcal{D}^*(u_j^i(\bm{\mathrm{x}},t;p)))h(\bm{\mathrm{x}},t)d\bm{\mathrm{x}}dt&=\int_{\Omega_a\times(0,T)}\mathcal{N}(\bm{\Theta}\cdot \mathcal{D}^*(u_j^i(\bm{\mathrm{x}},t;q)))h(\bm{\mathrm{x}},t)d\bm{\mathrm{x}}dt,\\
\end{split}
\end{eqnarray}
then $q=p$ in $\Omega$.\\

\textbf{Theorem 2.} Let $\{\phi_j\}_{j=1}^{\infty}\in C(\Omega)$ be a complete set in $L^2(\Omega)$, $v\in C^1(0,T)$ and $h\in C_0^{2,1}(\Omega_a\times(0,T))$ be given nonzero nonnegative functions, and $v$ satisfies $v(0)=0$. Assume $p(\bm{\mathrm{x}})$, $q(\bm{\mathrm{x}})\in C(\overline\Omega)$, $p$, $q\geq 0$ on $\Omega$. Let $u_j^i(\bm{\mathrm{x}},t;p)$, $u_j^i(\bm{\mathrm{x}},t;q)$ be the classical solutions of problem (\ref{mix2}) corresponding to the input sources $\phi_j(\bm{\mathrm{x}})v_i(t)$ ($i=1,2$; $j=1,2,\cdots$) with the reation coefficients $p$ and $q$ respectively. If we choose 
\begin{eqnarray*}
\begin{split}
v_1=v,\ \ \  v_2={_{0}D}_t^{\alpha}v+\sum\limits_{k=1}^{m}p_{k}({_{0}D}_t^{\alpha_k}v),
\end{split}
\end{eqnarray*}
such that
\begin{eqnarray}
\small
\label{condtion 1}
\begin{split}
\int_{\Omega_a\times(0,T)}\mathcal{N}(\bm{\Theta}\cdot \mathcal{D}^*(u_j^i(\bm{\mathrm{x}},t;p)))h(\bm{\mathrm{x}},t)d\bm{\mathrm{x}}dt&=\int_{\Omega_a\times(0,T)}\mathcal{N}(\bm{\Theta}\cdot \mathcal{D}^*(u_j^i(\bm{\mathrm{x}},t;q)))h(\bm{\mathrm{x}},t)d\bm{\mathrm{x}}dt,\\
\end{split}
\end{eqnarray}
then $q=p$ in $\Omega$. \\

\textbf{Remark 1.} Notice that if "=" hold in (\ref{kernel function}), i.e., the kernel function is given by
\begin{eqnarray}
\label{}
\begin{split}
\gamma(\bm{\mathrm{x}},\bm{\mathrm{y}})=\frac{\gamma_{*}}{\mid \bm{\mathrm{x}}-\bm{\mathrm{y}}\mid^{d+2\beta}},\ \ \ \text{for}\ \ \bm{\mathrm{y}}\in B_\varepsilon(\bm{\mathrm{x}}),
\end{split}
\end{eqnarray}
and we choose $\varepsilon=\infty$, then the operator $\mathcal{L}$ is simplified as fractional Laplacian $-(-\Delta)^\beta$ \cite{J. Jia2018}.

The paper is organized as follows. In Section 2, the preliminary is used to introduce the concept of nonlocal calculus. In Section 3, we prove that the uniqueness theorems for NDE and MTTFNDE.

\section{Preliminary}

In this section, we briefly review the concepts of nonlocal calculus that are useful in what follows. The principal goal is to develop a vector calculus for nonlocal operators that mimics the classical vector calculus for differential operators, refer to \cite{Q Du,Du2012}.

The action of the nonlocal divergence operator $\mathcal{D}(\bm{{\nu}}):\mathbb{R}^d\rightarrow \mathbb{R}$ on $\bm{\mathrm{\nu}}$ is defined as
\begin{eqnarray}
\label{Pre1}
\begin{split}
\mathcal{D}(\bm{\mathrm{\nu}})(\bm{\mathrm{x}}):=\int_{\mathbb{R}^d}(\bm{\mathrm{\nu}}(\bm{\mathrm{x}},\bm{\mathrm{y}})+\bm{\nu}(\bm{\mathrm{y}},\bm{\mathrm{x}}))\cdot \bm{\alpha}(\bm{\mathrm{x}},\bm{\mathrm{y}})d\bm{\mathrm{y}},\ \ \ \text{for}\ \ \bm{\mathrm{x}}\in \mathbb{R}^d,
\end{split}
\end{eqnarray}
where the vector mappings $\bm{\mathrm{\nu}}(\bm{\mathrm{x}},\bm{\mathrm{y}}),\ \bm{\alpha}(\bm{\mathrm{x}},\bm{\mathrm{y}}):\mathbb{R}^d\times \mathbb{R}^d\rightarrow \mathbb{R}^k$ with $\bm{\alpha}$ antisymmetric, i.e., $\bm{\alpha}(\bm{\mathrm{x}},\bm{\mathrm{y}})=-\bm{\alpha}(\bm{\mathrm{y}},\bm{\mathrm{x}})$.

Given the mapping $u(\bm{\mathrm{x}}):\mathbb{R}^d\rightarrow \mathbb{R}$, the adjoint operator $\mathcal{D}^*$ corresponding to $\mathcal{D}$ is the operator whose action on $u$ is given by
\begin{eqnarray}
\label{Pre2}
\begin{split}
\mathcal{D}^*(u)(\bm{\mathrm{x}},\bm{\mathrm{y}}):=-(u(\bm{\mathrm{y}})-u(\bm{\mathrm{x}}))\cdot \bm{\alpha}(\bm{\mathrm{x}},\bm{\mathrm{y}}),\ \ \ \text{for}\ \ \bm{\mathrm{x}},\bm{\mathrm{y}}\in \mathbb{R}^d,
\end{split}
\end{eqnarray}
 where $\mathcal{D}^*(u):\mathbb{R}^d\times \mathbb{R}^d\rightarrow \mathbb{R}^k$. In fact, $-\mathcal{D}^*$ denotes a nonlocal gradient.

We can see that if $\bm{\Theta}(\bm{\mathrm{x}},\bm{\mathrm{y}})=\bm{\Theta}(\bm{\mathrm{y}},\bm{\mathrm{x}})$ denotes a second-order symmetric definite tensor satisfying $\bm{\Theta}=\bm{\Theta}^{T}$, then
\begin{eqnarray}
\label{Pre3}
\begin{split}
\mathcal{D}(\bm{\Theta}\cdot\mathcal{D}^* u)(\bm{\mathrm{x}}):=-2\int_{\mathbb{R}^d}(u(\bm{\mathrm{y}})-u(\bm{\mathrm{x}}))\bm{\alpha}(\bm{\mathrm{x}},\bm{\mathrm{y}})\cdot(\bm{\Theta}\cdot\bm{\alpha}(\bm{\mathrm{x}},\bm{\mathrm{y}}))d\bm{\mathrm{y}},\ \ \ \text{for}\ \ \bm{\mathrm{x}}\in \mathbb{R}^d,
\end{split}
\end{eqnarray}
where $\mathcal{D}(\bm{\Theta}\cdot\mathcal{D}^* u):\mathbb{R}^d\rightarrow \mathbb{R}$. Comparing with $(\ref{mi1})$, we see that
\begin{eqnarray*}
\begin{split}
-\mathcal{L}u=\mathcal{D}(\bm{\Theta}\cdot\mathcal{D}^* u),\ \ \ \ \text{with}\ \ \gamma=\bm{\alpha}\cdot(\bm{\Theta}\cdot\bm{\alpha}).
\end{split}
\end{eqnarray*}

Given an open subset $\Omega\subset \mathbb{R}^d$, the corresponding interaction domain is defined by
\begin{eqnarray*}
\begin{split}
\Omega_{\mathcal{I}}:=\{\bm{\mathrm{y}}\in \mathbb{R}^d\setminus\Omega\ \ \ \text{such that}\ \ \alpha(\bm{\mathrm{x}},\bm{\mathrm{y}})\neq 0\ \ \text{for}\ \bm{\mathrm{x}} \in\Omega\}.
\end{split}
\end{eqnarray*}
So that $\Omega_{\mathcal{I}}$ consists of those points outside of $\Omega$ that interact with points in $\Omega$. Then, the corresponding to the divergence operator $\mathcal{D}(\bm{\nu}):\mathbb{R}^d\rightarrow \mathbb{R}$ defined in (\ref{Pre1}), we also define the action of the nonlocal interaction operator $\mathcal{N}(\bm{\nu}):\mathbb{R}^d\rightarrow \mathbb{R}$ on $\bm{\nu}$ by
\begin{eqnarray}
\label{Pre4}
\begin{split}
\mathcal{N}(\bm{\nu})(\bm{\mathrm{x}}):=-\int_{\Omega\cup\Omega_I}(\bm{\nu}(\bm{\mathrm{x}},\bm{\mathrm{y}})+\bm{\nu}(\bm{\mathrm{y}},\bm{\mathrm{x}}))\cdot\bm{\alpha}(\bm{\mathrm{x}},\bm{\mathrm{y}})d\bm{\mathrm{y}},\ \ \ \ \text{for}\ \bm{\mathrm{x}}\in\Omega_{\mathcal{I}}.
\end{split}
\end{eqnarray}
In \cite{Q Du}, it is shown that $\int_{\Omega_{\mathcal{I}}} \mathcal{N}(\bm{\nu})d\bm{\mathrm{x}}$ can ba viewed as a nonlocal flux out of $\Omega$ into $\Omega_{\mathcal{I}}$.

With $\mathcal{D}$ and $\mathcal{N}$ defined in (\ref{Pre1}) and (\ref{Pre4}), respectively, we have the nonlocal Gauss theorem
\begin{eqnarray}
\label{Pre5}
\begin{split}
\int_{\Omega}\mathcal{D}(\bm{\nu})d\bm{\mathrm{x}}=\int_{\Omega_{\mathcal{I}}}\mathcal{N}(\bm{\nu})d\bm{\mathrm{x}}.
\end{split}
\end{eqnarray}

Next, let $u(\bm{\mathrm{x}})$ and $\bm{\nu}(\bm{\mathrm{x}})$ denote scalar functions. Then we can show that the nonlocal divergence theorem (\ref{Pre5}) implies the nonlocal Green's first identity
\begin{eqnarray}
\label{Pre6}
\begin{split}
\int_{\Omega}\bm{\nu}\mathcal{D}(\bm{\Theta}\cdot\mathcal{D}^*u)d\bm{\mathrm{x}}-\int_{\Omega\cup\Omega_{\mathcal{I}}}\int_{\Omega\cup\Omega_{\mathcal{I}}}(\mathcal{D}^*\bm{\nu})\cdot(\bm{\Theta}\cdot\mathcal{D}^*u)d\bm{\mathrm{y}}d\bm{\mathrm{x}}
=\int_{\Omega_{\mathcal{I}}}\bm{\nu}\mathcal{N}(\bm{\Theta}\cdot\mathcal{D}^*u)d\bm{\mathrm{x}}.
\end{split}
\end{eqnarray}


\section{The uniqueness of the IRCP for NDE and MTTFNDE}
In this section, we show that the measurement data $(\ref{me1})$ can determine the reaction coefficient $q$ uniquely for NDE and MTTFNDE. In order to prove the uniqueness, the nonlocal maximum principle will be established here (see also \cite{Bucur2016} for fractional Laplacian case). \\

\textbf{Lemma 1.}\ (Weak Maximum Principle) Assume $u\in C^{2,1}_{0}(\Omega\times(0,T))$, if $\frac{\partial u}{\partial t}-\mathcal{L}u+q(\bm{\mathrm{x}})u \geq 0$ in $\Omega\times(0,T)$, and $u\geq 0$ in $(\mathbb{R}^d\setminus \Omega)\times(0,T)$, then we have $u\geq 0$ in $\Omega\times(0,T)$.\\

\textbf{Proof}. Assume now  by contradiction that the minimal point $(\bm{\mathrm{x_0}},t_0)\in \Omega\times(0,T)$ is attained and satisfies $u(\bm{\mathrm{x_0}},t_0)<0$, since $u$ is nonnegative outside $\Omega\times(0,T)$. Then $u(\bm{\mathrm{x_0}},t_0)$ is a minimum in $\mathbb{R}^d\times(0,T)$ and deduces that $\frac{\partial u}{\partial t}\mid_{(\bm{\mathrm{x_0}},t_0)}=0$. We set $r=dist(\bm{\mathrm{x_0}},\partial\Omega)$ and $B_r(\bm{\mathrm{x_0}})$ denotes the center of the circle is $\bm{\mathrm{x_0}}$, with a radius of $r$. Due to $u(\bm{\mathrm{x_0}},t_0)$ is a minimum, we have $u(\bm{\mathrm{x_0}},t_0)-u(\bm{\mathrm{y}},t_0)\leq 0$, for $\bm{\mathrm{y}}\in B_{2r}(\bm{\mathrm{x_0}})$. If $\bm{\mathrm{y}}\in \mathbb{R}^d\setminus B_{2r}(\bm{\mathrm{x_0}})$, then $\mid \bm{\mathrm{x}}-\bm{\mathrm{y}}\mid \geq\mid \bm{\mathrm{y}}-\bm{\mathrm{x_0}}\mid-\mid \bm{\mathrm{x}}-\bm{\mathrm{x_0}}\mid\geq r$ and $u(\bm{\mathrm{y}},t_0)\geq 0$. Thus, according to $(\ref{kernel function})$,
\begin{eqnarray*}
\begin{split}
0&\leq(\frac{\partial u}{\partial t}-\mathcal{L}u+q(\bm{\mathrm{x}})u)\mid_{(\bm{\mathrm{x_0}},t_0)}\\
&=\frac{\partial u}{\partial t}\mid_{(\bm{\mathrm{x_0}},t_0)}+2\int_{\mathbb{R}^d}(u(\bm{\mathrm{x_0}},t_0)-u(\bm{\mathrm{y}},t_0)) \gamma(\bm{\mathrm{x_0}},\bm{\mathrm{y}})d\bm{\mathrm{y}}+q(\bm{\mathrm{x_0}})u(\bm{\mathrm{x_0}},t_0)\\
&\leq 2\int_{\mathbb{R}^d}\frac{\gamma_*(u(\bm{\mathrm{x_0}},t_0)-u(\bm{\mathrm{y}},t_0))}{\mid \bm{\mathrm{x_0}}-\bm{\mathrm{y}}\mid^{n+2\beta}}d\bm{\mathrm{y}}\\
&=2\gamma_*\int_{B_{2r}(\bm{\mathrm{x_0}})}\frac{u(\bm{\mathrm{x_0}},t_0)-u(\bm{\mathrm{y}},t_0)}{\mid \bm{\mathrm{x_0}}-\bm{\mathrm{y}}\mid^{n+2\beta}}d\bm{\mathrm{y}}+2\gamma_*\int_{\mathbb{R}^d\setminus B_{2r}(\bm{\mathrm{x_0}})}\frac{u(\bm{\mathrm{x_0}},t_0)-u(\bm{\mathrm{y}},t_0)}{\mid \bm{\mathrm{x_0}}-\bm{\mathrm{y}}\mid^{n+2\beta}}d\bm{\mathrm{y}}\\
&\leq 2\gamma_*\int_{\mathbb{R}^d \setminus B_{2r}(\bm{\mathrm{x_0}})}\frac{u(\bm{\mathrm{x_0}},t_0)}{\mid \bm{\mathrm{x_0}}-\bm{\mathrm{y}}\mid^{n+2\beta}}d\bm{\mathrm{y}}<0.
\end{split}
\end{eqnarray*}
It leads to contradictions, so we can get $u\geq 0$ in $\Omega\times(0,T)$.\\

\textbf{Lemma 2.}\ (Strong Maximum Principle) Assume $u\in C^{2,1}_{0}(\Omega\times(0,T))$, if $\frac{\partial u}{\partial t}-\mathcal{L}u+q(\bm{\mathrm{x}})u \geq 0$ in $\Omega\times(0,T)$, and $u\geq0$ in $(\mathbb{R}^d\setminus \Omega)\times(0,T)$, then  $u\geq0$ in $\Omega\times(0,T)$, unless u vanishes identically.\\

\textbf{Proof}. We observe that we already know that $u\geq0$ in $\mathbb{R}^d\times(0,T)$ according to the Lemma 1. Hence, if $u$ is not strictly positive, there exists $(\bm{\mathrm{x_0}},t_0)\in \Omega\times(0,T)$ such that $u(\bm{\mathrm{x_0}},t_0)=0$. This gives that
\begin{eqnarray*}
\begin{split}
0&\leq(\frac{\partial u}{\partial t}-\mathcal{L}u+q(\bm{\mathrm{x}})u)\mid_{(\bm{\mathrm{x_0}},t_0)}\\
&=\frac{\partial u}{\partial t}\mid_{(\bm{\mathrm{x_0}},t_0)}+2\int_{\mathbb{R}^d}(u(\bm{\mathrm{x_0}},t_0)-u(\bm{\mathrm{y}},t_0)) \gamma(\bm{\mathrm{x_0}},\bm{\mathrm{y}})d\bm{\mathrm{y}}+q(\bm{\mathrm{x_0}})u(\bm{\mathrm{x_0}},t_0)\\
&=-2\int_{\mathbb{R}^d}u(\bm{\mathrm{y}},t_0)\gamma(\bm{\mathrm{x_0}},\bm{\mathrm{y}})d\bm{\mathrm{y}}\leq 0,\\
\end{split}
\end{eqnarray*}
then we can get $u\equiv 0$ in $\Omega\times(0,T)$, the conclusion is established.\\

\textbf{Proof of Theorem 1}. Notice that $h\in C_0^{2,1}(\Omega_a\times(0,T))$ be a given nonzero nonnegative function. Then we set $h_0=h$ on $(\Omega_a\times(0,T))$ and $h_0=0$ on $(\mathbb{R}^d\setminus(\Omega\cup\Omega_a))\times(0,T)$, and introduce the function $ w(x,t;q)$ as the solution of the following adjoint problem
\begin{eqnarray}
\label{adjoint problem}
\begin{cases}
\begin{split}
-\frac{\partial w}{\partial t}-\mathcal{L}w+q(\bm{\mathrm{x}})w&=0,\ \ \ \text{in}\ \ \Omega\times(0, T),\\
w&=h_0(\bm{\mathrm{x}},t),\ \ \ \text{on}\ \ (\mathbb{R}^d\setminus\Omega)\times(0, T),\\
w&=0,\ \ \ \text{in}\ \ \Omega\times\{T\}.\\
\end{split}
\end{cases}
\end{eqnarray}
In fact, by using the transform formula $\tilde{w}(\bm{\mathrm{x}},t)=w(\bm{\mathrm{x}},T-t)$
\begin{eqnarray}
\label{adjoint problem con}
\begin{cases}
\begin{split}
\frac{\partial \tilde{w}}{\partial t}-\mathcal{L}\tilde{w}+q(\bm{\mathrm{x}})\tilde{w}&=0,\ \ \ \text{in}\ \ \Omega\times(0, T),\\
\tilde{w}&=h_0(\bm{\mathrm{x}},T-t),\ \ \ \text{on}\ \  (\mathbb{R}^d\setminus\Omega)\times(0, T),\\
\tilde{w}&=0,\ \ \ \text{in}\ \ \Omega\times\{0\}.\\
\end{split}
\end{cases}
\end{eqnarray}
Since we choose $v_1=v$, from  $(\ref{mix1})$, $(\ref{adjoint problem})$ and the nonlocal Green's formula $(\ref{Pre6})$ we compute
\begin{eqnarray*}\cdot
\small
\label{Qresidual-eq}
\begin{split}
&\int_{0}^T \int_{\Omega}\phi_{j}(\bm{\mathrm{x}})v(t)w(\bm{\mathrm{x}},t;q)d\bm{\mathrm{x}}dt\\= & \int_{0}^T\int_{\Omega}[\frac{\partial u_j^1(\bm{\mathrm{x}},t;q) }{\partial t}-\mathcal{L}u_j^1(\bm{\mathrm{x}},t;q)+q(\bm{\mathrm{x}})u_j^1(\bm{\mathrm{x}},t;q)]w(\bm{\mathrm{x}},t;q)d\bm{\mathrm{x}}dt\\
=&-\int_0^T \int_{\Omega} (\frac{\partial w}{\partial t})u_j^1d\bm{\mathrm{x}}dt-\int_0^T\int_{\Omega}\mathcal{L}wu_j^1dxdt+\int_0^T \int_{\Omega}q(\bm{\mathrm{x}})u_j^1 wd\bm{\mathrm{x}}dt\\&+\int_0^T\int_{\mathbb{R}^d\setminus\Omega}w\mathcal{N}(\bm{\Theta}\cdot \mathcal{D}^*(u_j^1(\bm{\mathrm{x}},t;q)))d\bm{\mathrm{x}}dt-\int_0^T\int_{\mathbb{R}^d\setminus\Omega}u_j^1\mathcal{N}(\bm{\Theta}\cdot \mathcal{D}^*(w(\bm{\mathrm{x}},t;q)))d\bm{\mathrm{x}}dt\\
=&\int_{0}^T\int_{\Omega}[-\frac{\partial w}{\partial t}-\mathcal{L}w+q(\bm{\mathrm{x}})w]u_j^1d\bm{\mathrm{x}}dt+\int_0^T\int_{\Omega_a} h(\bm{\mathrm{x}},t) \mathcal{N}(\bm{\Theta}\cdot \mathcal{D}^*(u_j^1(\bm{\mathrm{x}},t;q)))d\bm{\mathrm{x}}dt\\
=&\int_0^T\int_{\Omega_a} h(\bm{\mathrm{x}},t)\mathcal{N}(\bm{\Theta}\cdot \mathcal{D}^*(u_j^1(\bm{\mathrm{x}},t;q)))d\bm{\mathrm{x}}dt.\\
\end{split}
\end{eqnarray*}
 Similarly, set $v_2(t)=v'(t)$, it follows that
\begin{eqnarray*}
\begin{split}
\int_{0}^T\int_{\Omega} \phi_{j}(\bm{\mathrm{x}})v'(t)w(\bm{\mathrm{x}},t;q)d\bm{\mathrm{x}}dt=\int_0^T\int_{\Omega_a} \mathcal{N}(\bm{\Theta}\cdot \mathcal{D}^*(u_j^2(\bm{\mathrm{x}},t;q)))h(\bm{\mathrm{x}},t)d\bm{\mathrm{x}}dt.
\end{split}
\end{eqnarray*}
For $p(\bm{\mathrm{x}})$ and the corresponding function $w(\bm{\mathrm{x}},t;p)$ given by $(\ref{adjoint problem})$ ,we see that
\begin{eqnarray*}
\begin{split}
\int_{0}^T\int_{\Omega} \phi_{j}(\bm{\mathrm{x}})v(t)w(\bm{\mathrm{x}},t;p)d\bm{\mathrm{x}}dt&=\int_0^T\int_{\Omega_a}\mathcal{N}(\bm{\Theta}\cdot \mathcal{D}^*(u_j^1(\bm{\mathrm{x}},t;p)))h(\bm{\mathrm{x}},t)d\bm{\mathrm{x}}dt;\\
\int_{0}^T \int_{\Omega} \phi_{j}(\bm{\mathrm{x}})v'(t))w(\bm{\mathrm{x}},t;p)d\bm{\mathrm{x}}dt&=\int_0^T\int_{\Omega_a}\mathcal{N}(\bm{\Theta}\cdot \mathcal{D}^*(u_j^2(\bm{\mathrm{x}},t;p)))h(\bm{\mathrm{x}},t)d\bm{\mathrm{x}}dt.
\end{split}
\end{eqnarray*}
Then we see from  $(\ref{condtion})$ that
\begin{eqnarray*}
\begin{split}
\int_{0}^T\int_{\Omega} \phi_{j}(\bm{\mathrm{x}})v(t)w(\bm{\mathrm{x}},t;q)d\bm{\mathrm{x}}dt=&\int_{0}^T \int_{\Omega} \phi_{j}(\bm{\mathrm{x}})v(t)w(\bm{\mathrm{x}},t;p)d\bm{\mathrm{x}}dt;\\
\int_{0}^T\int_{\Omega} \phi_{j}(\bm{\mathrm{x}})v'(t)w(\bm{\mathrm{x}},t;q)d\bm{\mathrm{x}}dt=&\int_{0}^T\int_{\Omega} \phi_{j}(\bm{\mathrm{x}})v'(t)w(\bm{\mathrm{x}},t;p)d\bm{\mathrm{x}}dt.
\end{split}
\end{eqnarray*}
By the completeness of $\{\phi_j(\bm{\mathrm{x}})\}_{j=1}^{\infty}$, we obtain
\begin{eqnarray}
\label{eq11}
\begin{split}
\int_{0}^T v(t)w(\bm{\mathrm{x}},t;q)dt=&\int_{0}^T v(t)w(\bm{\mathrm{x}},t;p)dt;\\
\int_{0}^T v'(t)w(\bm{\mathrm{x}},t;q)dt=&\int_{0}^T v'(t)w(\bm{\mathrm{x}},t;p)dt.\\
\end{split}
\end{eqnarray}
Multiplying equation (\ref{adjoint problem}) by $v$, integrating by parts over $(0,T)$, we find
\begin{eqnarray*}
\begin{split}
-\int_{0}^T\frac{\partial w(\bm{\mathrm{x}},t;p)}{\partial t}v(t)dt+\int_0^T \mathcal{D}(\bm{\Theta}\cdot \mathcal{D}^*w(\bm{\mathrm{x}},t;p))v(t)dt+\int_0^T p(\bm{\mathrm{x}})w(\bm{\mathrm{x}},t;p)v(t)dt=0;\\
-\int_{0}^T\frac{\partial w(\bm{\mathrm{x}},t;q)}{\partial t}v(t)dt+\int_0^T \mathcal{D}(\bm{\Theta}\cdot \mathcal{D}^*w(\bm{\mathrm{x}},t;q))v(t)dt+\int_0^T q(\bm{\mathrm{x}})w(\bm{\mathrm{x}},t;q)v(t)dt=0;\\
\end{split}
\end{eqnarray*}
Then deduce that
\begin{eqnarray}
\label{eq2}
\begin{split}
\int_{0}^T w(\bm{\mathrm{x}},t;p)v'(t)dt+ \mathcal{D}(\bm{\Theta}\cdot \mathcal{D}^*(\int_0^Tw(\bm{\mathrm{x}},t;p)v(t)dt))+\int_0^T p(\bm{\mathrm{x}})w(\bm{\mathrm{x}},t;p)v(t)dt=0;\\
\int_{0}^T w(\bm{\mathrm{x}},t;q)v'(t)dt+ \mathcal{D}(\bm{\Theta}\cdot \mathcal{D}^*(\int_0^Tw(\bm{\mathrm{x}},t;q)v(t)dt))+\int_0^T q(\bm{\mathrm{x}})w(\bm{\mathrm{x}},t;q)v(t)dt=0;\\
\end{split}
\end{eqnarray}
The two expressions of $(\ref{eq2})$ are subtracted from each other, and using (\ref{eq11}) we have
\begin{eqnarray}
\label{eq3}
\begin{split}
(p-q)\int_0^T w(\bm{\mathrm{x}},t;q)v(t)dt=0,
\end{split}
\end{eqnarray}
The strong maximum principle of the Lemma 2 can be applied to deduce that $w(\bm{\mathrm{x}},t;q)>0$, then $q=p$ in $\Omega$.\\

Next, similar to the NDE case, we also prove the corresponding maximum principle and establish the uniqueness of IRCP for MTTFNDE.\\

\textbf{Lemma 3.}\ (\cite{X. Yang2018} Lemma 1) Assume that $f\in C[0,T]\cap C^1[0,T]$, attain its minimum over the interval [0,T] at a point to $t_0\in (0,T]$, then
\begin{eqnarray}
\label{eq8}
\begin{split}
_{0}D_t^{\alpha}f(t_0)\leq0.
\end{split}
\end{eqnarray}

By using Lemma 3, the nonlocal weak maximum principle and strong maximum principle can be obtained similarly for MTTFNDE.\\

\textbf{Lemma 4.}\ (Weak Maximum Principle) Assume $u\in C^{2,1}_{0}(\Omega\times(0,T))$, if
\begin{eqnarray}
\begin{split}
{_{0}D}_t^{\alpha}u+\sum_{i=1}^{m}p_k(_{0}D_t^{\alpha_k}u)-\mathcal{L}u+q(\bm{\mathrm{x}})u \geq 0,\ \ \  \mathrm{in}\ \Omega\times(0,T),
\end{split}
\end{eqnarray}
and $u\geq 0$ in $(\mathbb{R}^d\setminus \Omega)\times(0,T)$, then we can deduce that $u\geq 0$ in $\Omega\times(0,T)$.\\

\textbf{Lemma 5.}\ (Strong Maximum Principle) Assume $u\in C^{2,1}_{0}(\Omega\times(0,T))$, if
\begin{eqnarray}
\begin{split}
{_{0}D}_t^{\alpha}u+\sum_{i=1}^{m}p_k(_{0}D_t^{\alpha_k})u-\mathcal{L}u+q(x)u \geq 0, \ \ \mathrm{in}\ \  \Omega\times(0,T),
\end{split}
\end{eqnarray}
and $u\geq0$ in $(\mathbb{R}^d\setminus \Omega)\times(0,T)$, then we have $u\geq0$ in $\Omega\times(0,T)$, unless $u$ vanishes identically.\\

\textbf{Proof of Theorem 2}. Since the proof of Theorem 2 is very similar to the one in Theorem 1 by applying Maximum Principle, we omit it here.

\section*{Acknowledgments}
We acknowledge the support of the NSF of China (11301168).

\label{ssec:Conclusions}

\end{document}